\documentclass[11pt]{article}
\usepackage{mathrsfs}
\usepackage{latexsym,lineno}
\usepackage{epsfig}
\usepackage{color}
\usepackage{amsmath}\usepackage{verbatim}\usepackage{epsf}
\usepackage{amsthm}\usepackage{graphicx, float}\usepackage{graphicx}
\usepackage{amsfonts}\usepackage{amssymb}\usepackage{graphpap}
\usepackage{epic}\usepackage{curves}
\usepackage{tikz}
\usepackage{subfigure}

\newtheorem{fact}{Theorem}
\newtheorem{claim}{Claim}

\newtheorem{theorem}[fact]{Theorem}

\newtheorem{definition}{Definition}

\newtheorem{lemma}{Lemma}

\title{\bf  Redundant edges in Ramsey graphs\thanks{Supported in part by NSFC and Shanghai Sailing Program(19YF1435500).} }
\author{Ye Wang$^{a}$, \hspace{2mm} Yusheng Li$^{b}$, \hspace{2mm} Yan Li$^{b}$\footnote{Corresponding author: 1610521@tongji.edu.cn} \\
{\small $^{a}$ School of Stats \& Math, \; Shanghai Lixin University of Accounting and Finance} \\
{\small Shanghai 201209, China}\\
{\small $^{b}$ School of Mathematical Sciences, Tongji University,\; Shanghai 200092, China} \\
{\small Email: wangye@sfu.edu.cn,\;li\_yusheng@tongji.edu.cn,\;1610521@tongji.edu.cn }   }
\date{}
\begin{document}
\maketitle
\begin{abstract}

For graphs $G$, $F$ and $H$, let $G\rightarrow (F,H)$ signify
that any edge coloring of $G$ in red and blue contains a red $F$ or a blue $H$.
The Ramsey number $R(F,H)=\min\{r|\; K_r\rightarrow (F,H)\}$.
In this note, we consider redundant edges in Ramsey graphs,
which are associate with critical Ramsey numbers.
For an integer $k\ge 1$, let ${\mathbb G}=\{G_k,G_{k+1},\dots \}$ be a class of graphs with $\delta(G_n)\ge 1$.
We define the critical Ramsey number $R_{\mathbb G}(F,H)$ with respect to $\mathbb G$ to be
$\max\{n|\; K_r\setminus G_n \rightarrow(F,H),\,G_n\in{\mathbb G}\big\}$, where $r=R(F,H)$.
We shall determine some $R_{\mathbb G}(F,H)$, where ${\mathbb G}$ consists of stars, matchings and complete graphs, respectively.
\medskip


\end{abstract}

\section{Introduction}
\indent

For graphs $G$, $F$ and $H$, let $G\rightarrow (F,H)$ signify
that any edge coloring of $G$ in red and blue contains a red $F$ or a blue $H$.
Then $G \not\rightarrow (F,H)$ means that there is a red-blue edge coloring of $G$
that contains neither red $F$ nor blue $H$. The Ramsey number $R(F,H)=\min\{r|\; K_r\rightarrow (F,H)\}$.
If there is a red-blue edge coloring of $G$ such that there is neither red $F$ nor blue $H$,
we call such a coloring to be an $(F,H)$-free coloring and the graph $G$ is called an $(F,H)$-free graph.
In particular, for $r=R(F,H)$, an $(F,H)$-free coloring of $K_{r-1}$ is called a Ramsey coloring of $K_{r-1}$,
and the corresponding monochromatic graphs are called critical graphs, where the two  monochromatic graphs are complementary.

For vertex disjoint graphs $F$ and $H$, denote by $F+H$ the join of $F$ and $H$
obtained by additional edges connecting $V(F)$ and $V(H)$ completely, and $F\cup H$ the disjoint union of $F$ and $H$.
If $V(F)\subseteq V(K_n)$,
let $K_n\setminus F$ be the graph obtained from $K_n$ by deleting the edges of $F$ from $K_n$.

In this note, we consider redundant edges in Ramsey graphs,
which are associate with critical Ramsey numbers.
Let $S_n=K_{1,n}$ and $M_n=nK_2$ be a star and a matching of $n$ edges, respectively.
Let $P_n$ be path on $n$ vertices. Denote
\[
\mathbb S=\{S_1,S_2,\dots\},\;\mathbb M=\{M_1,M_2,\dots\},\;\mathbb P=\{P_2,P_3,\dots\},\;\mathbb K=\{K_2,K_3,\dots\}.
\]
Here the indices of graphs start at different $k$ to coincide with conventional notation in graph theory. The indices of stars and matchings start at $k=1$,
and those of paths and complete graphs start at $k=2$.
In the following definition, we denote by $\delta(G)$ the minimum degree of $G$.
One can also find this definition in our recent paper \cite{wang}.
\begin{definition}\label{d2}
Let ${\mathbb G}=\{G_k,G_{k+1},\dots \}$ be a class of graphs,
where the fixed integer $k\ge 1$ and $\delta(G_n)\ge 1$ for $G_n\in \mathbb G$.
Define the critical Ramsey number $R_{\mathbb G}(F,H)$ with respect to $\mathbb G$ as
\[
R_{\mathbb G}(F,H)=\max\big\{n|\; K_r\setminus G_n \rightarrow(F,H),\,G_n\in{\mathbb G}\big\},
\]
where $r=R(F,H)$.
\end{definition}
We shall call $R_{\mathbb S}(F,H)$ the star-critical Ramsey number,
$R_{\mathbb M}(F,H)$ the matching-critical Ramsey number, $R_{\mathbb P}(F,H)$ the path-critical Ramsey number
and $R_{\mathbb K}(F,H)$ the complete-critical Ramsey number, respectively.
Note that $n\ge 1$ for any $G_n\in\mathbb{G}$, thus we write $R_{\mathbb G}(F,H)=0$ if $K_r\setminus G_n\not\rightarrow(F,H)$ for any $G_n\in {\mathbb G}$.

Usually, we shall consider such class ${\mathbb G}$ that $G_n$ is a proper subgraph of $G_{n+1}$ for $n\geq k$.
Note that $R_{\mathbb G}(F,H)$ is related to the indices of graphs in $\mathbb G$. For example, if $K_2$ is the graph with maximum number of edges such that $K_r\setminus K_2\rightarrow(F,H)$,
then $R_{\mathbb S}(F,H)=R_{\mathbb M}(F,H)=1$ and $R_{\mathbb P}(F,H)=R_{\mathbb K}(F,H)=2$.

The concept of star-critical Ramsey number was introduced by Hook and Isaak \cite{hook} in a different way.
Let $K_n\sqcup S_k$ be the graph obtained from $K_n$ and an additional vertex $v$ by joining $v$ to $k$ vertices of $K_n$.
Hook and Isaak \cite{hook} defined the star-critical Ramsey number $r_*(F,H)$ as the smallest $k$
such that $K_{r-1}\sqcup S_k\rightarrow (F,H)$, where $r=R(F,H)$. Hook and Isaak \cite{hook} determined some star-critical Ramsey numbers,
and Zhang, Broersma and Chen \cite{zhang}, Hao and Lin \cite{lin} showed general lower bounds for $r_*(F,H)$.

From the definitions, we have an easy equality as
\begin{equation}\label{e1}
R_{\mathbb S}(F,H)=R(F,H)-1-r_*(F,H).
\end{equation}
Hence we have the following results from the known star-critical Ramsey numbers,
in which $T_n$ is a tree on $n$ vertices and $F_n=K_1+nK_2$.
\bigskip

\begin{center}
\begin{tabular}{|l|l|l|}   \hline
$r_*(F,H)$ &$R_{\mathbb S}(F,H)$& From\\ \hline
$r_*(T_n,K_m)=(n-1)(m-2)+1$, $n\geq 2,m\geq 2$ &$R_{\mathbb S}(T_n,K_m)=n-2$& \cite{hook} \\ \hline
$r_*(nK_2,mK_2)=m$, $n\geq m\geq 1$&$R_{\mathbb S}(nK_2,mK_2)=2n-2$& \cite{hook} \\ \hline
$r_*(nK_3,mK_3)=3n+2m-1$, $n\geq m\geq 1,n\geq 2$&$R_{\mathbb S}(nK_3,mK_3)=0$& \cite{hook} \\ \hline
$r_*(P_n,C_4)=3$, $n\ge 3$&$R_{\mathbb S}(P_n,C_4)=n-3$& \cite{hook} \\ \hline
$r_*(C_n,C_4)=5$, $n\geq 6$&$R_{\mathbb S}(C_n,C_4)=n-5$& \cite{wu}\\ \hline
$r_*(P_n,P_m)=\lceil m/2 \rceil$, $n\geq m\geq 4$
&$R_{\mathbb S}(P_n,P_m)=\left\{\begin{array}{cl} n-3, & \mbox{odd $m$} \\
   n-2, & \mbox{even $m$} \end{array} \right. $
& \cite{hook15}\\\hline
$r_*(nK_4,mK_3)=\left\{\begin{array}{cl} 4n+2m, & n\ge\max\{2,m\} \\
   3n+3m, & 2\le n< m \end{array} \right. $
&$R_{\mathbb S}(nK_4,mK_3)=0$
& \cite{li} \\\hline
$r_*(K_n,mK_2)=n+2m-3$, $m\geq 1,n\geq 3$ &$R_{\mathbb S}(K_n,mK_2)=0$& \cite{li}\\ \hline
$r_*(F_n,K_3)=2n+2$, $n\geq 2$ &$R_{\mathbb S}(F_n,K_3)=2n-2$& \cite{li}\\ \hline
$r_*(F_n,K_4)=4n+2$, $n\geq 4$ &$R_{\mathbb S}(F_n,K_4)=2n-2$& \cite{haghi}\\ \hline
\end{tabular}

\vspace{3mm}
 \bf{Table 1:} The known star-critical Ramsey numbers
\end{center}

Denote by $v(F)$ the order of $F$, $e(F)$ the size of $F$,
and $s(F)$ the chromatic surplus of $F$ which is the size of the smallest color class
in a $\chi (F)$-coloring of $F$.
We describe the connected graph $H$ with $v(H)\geq s(F)$ as $F$-good
if $R(F,H)=(\chi (F)-1)(v(H)-1)+s(F)$.
Let $\tau(F)$ be the minimum degree of some vertex of $U_1$ in $U_i$ for $2\le i\le \chi(F)$ under all vertex colorings of $F$ with $|U_1|=s(F)$
and the other color classes $U_2,\ldots,U_{\chi(F)}$.
Hence we have a general upper bound for $R_{\mathbb S}(F,H)$ by the lower bound for $r_*(F,H)$ in \cite{lin} and (\ref{e1}).

\begin{lemma}\cite{lin}\label{l4}
For graph $F$ with $\chi(F)\geq 2$ and connected graph $H$ of order $n\ge s(F)$,
if $H$ is $F$-good, then
\[
R_{\mathbb S}(F,H)\le \max\{s(F)-2,n+s(F)-\delta(H)-\tau(F)-1\}.
\]
\end{lemma}

Erd\H{o}s, Faudree, Rousseau and Schelp \cite{Erdos1978} introduced size Ramsey number to be
\[
\hat{r}(F,H)=\min\{e(G): \; G \rightarrow(F,H)\}.
\]

As usual, $\binom{m}{2}=0$ if the integer $m<2$.
In particular, $\binom{R_{\mathbb{K}}(F,H)}{2}=0$ if $R_{\mathbb{K}}(F,H)=0$, which means that $K_r\setminus K_n\not\to (F,H)$ for any $K_n\in \mathbb K$ as mentioned.
We shall omit the easy proof for the following result.
\begin{lemma}\label{l7}
For graphs $F$ and $H$ with $r=R(F,H)$,
\[
\hat{r}(F,H)\leq \binom{r}{2}-\binom{R_{\mathbb{K}}(F,H)}{2}.
\]
\end{lemma}

The following result is due to Chv\'atal, reported in Erd\H{o}s, Faudree, Rousseau and Schelp \cite{Erdos1978}.

\begin{lemma} \cite{Erdos1978}\label{l8}
For integers $m$ and $n$, $$\hat{r}(K_m,K_n)=\binom{r}{2}.$$
\end{lemma}
Hence from Lemma \ref{l7} and Lemma \ref{l8}, we have $R_{\mathbb K}(K_m,K_n)=0$.

The authors \cite{zhang} called $(F,H)$ Ramsey-full if $K_r\rightarrow (F,H)$, but $K_r\setminus K_2\not\rightarrow (F,H)$,
thus we call $(F,H)$ Ramsey-full on ${\mathbb G}$ if $R_{\mathbb G}(F,H)=0$.
From the known Ramsey-full pairs summarized in \cite{zhang}, we have that if $\mathbb G\in \{\mathbb S, \mathbb M, \mathbb P, \mathbb K\}$, then

\begin{equation*}
\begin{split}
&R_{\mathbb G}(K_m,K_n)=0 \quad \mbox{for any $m,n\ge 1$,}\\
&R_{\mathbb G}(K_n,mK_2)=0 \quad \mbox{for any $m\geq 1$ and $n\geq 3$,}\\
&R_{\mathbb G}(nK_3,mK_3)=0 \quad \mbox{for any $n\geq m\geq 1$ and $n\geq 2$,}\\
&R_{\mathbb G}(nK_4,mK_3)=0 \quad \mbox{for any $m\geq 1$ and $n\geq 2$}.
\end{split}
\end{equation*}

A path in a graph is called suspended if the degree of each internal vertex is two.
Bondy and Erd\H{o}s \cite{erdos73} proved that a long cycle $C_n$ (hence a long path) is $C_m$-good and $K_r(t)$-good.
Furthermore, Burr \cite{burr81} showed that $H$ is $F$-good for any fixed $F$ if $H$ contains a sufficiently long suspended path.
Note that $v(H)$ is the order of $H$ in following results.

\begin{lemma}\cite{burr81}\label{l1}
For graph $F$ and connected graph $H$,
let $H_n$ be a graph containing a suspended path of length $n-v(H)+1$ obtained from $H$
putting $n-v(H)$ extra vertices to an edge.
For sufficiently large $n$, $H_n$ is $F$-good.
\end{lemma}

\begin{theorem}\label{t1}
For graph $F$ with $\chi(F)\geq 2$ and connected graph $H$,
let $H_n$ be a graph containing a suspended path of length $n-v(H)+1$ obtained from $H$
putting $n-v(H)$ extra vertices to an edge. Then
$$
R_{\mathbb S}(F,H_n)=n+ C_n(F,H),
$$
where $|C_n(F,H)|\le v^2(F)+v(H)$.
\end{theorem}

From the value $r_*(P_n,P_m)$ determined in \cite{hook}, we have $R_{\mathbb S}(P_n,P_m)=n-3$ if $m$ is odd,
and $R_{\mathbb S}(P_n,P_m)=n-2$ otherwise.
We do not expect that $C_n(F,H)$ is an absolute constant
as the values of $R_{\mathbb S}(P_n,P_m)$.

\begin{theorem} \label{t2}
Let $m$ and $n$ be positive integers. Then
\[
R_{\mathbb S}(S_m,S_n)= \left\{  \begin{array}{cl}
0  & \mbox{if $m$ and $n$ are both even,}\\
m+n-2 & \mbox{otherwise}.
\end{array}      \right.
\]
\end{theorem}

\begin{theorem}\label{t3}
For positive integers $m$ and $n$,
\[
R_{\mathbb M}(S_m,S_n)= \left\{  \begin{array}{cl}
0  & \mbox{if $m$ and $n$ are both even,}\\
\Big\lfloor\frac{m+n-1}{2}\Big\rfloor & \mbox{otherwise}.
\end{array}      \right.
\]
\end{theorem}
\begin{theorem}\label{t4}
Let $m$ and $n$ be integers with $n\geq m\geq 1$ and $n\geq 2$. Then
\[
R_{\mathbb M}(M_m,M_n)=\Big\lfloor\frac{2n+m-1}{2}\Big\rfloor.
\]
\end{theorem}

In addition to $R_{\mathbb K}(S_1,S_1)=0$, we shall determine all other $R_{\mathbb K}(S_m,S_n)$.

\begin{theorem}\label{t5}
Let $m$ and $n$ be positive integers with $m+n\ge 3$. Then
\[
R_{\mathbb K}(S_m,S_n)= \left\{  \begin{array}{cl}
0  & \mbox{if $m$ and $n$ are both even,}\\
m+n-1 & \mbox{otherwise}.
\end{array}      \right.
\]
\end{theorem}
\begin{theorem}\label{t6}
Let $m$ and $n$ be integers with $n\geq m\geq 1$ and $n\geq 2$. Then
\[
R_{\mathbb K}(M_m,M_n)=n.
\]
\end{theorem}

Before proceeding to proofs, we need some notations.
For a red-blue edge colored $G$ in proofs,
the subgraph of $G$ induced by red edges is denoted by $G^R$. Similarly, the subgraph of $G$ induced by blue edges is denoted by $G^B$.
For a vertex $x$ of $G$, we denote the set of all neighbors of $x$ in $G^R$ by $N^R_G(x)$ and in $G^B$ by $N^B_G(x)$.
Let $d_G(x)$ (or $d(x)$ simply) denote the degree of vertex $x$ in $G$,
and $d^R_G(x)$, $d^B_G(x)$ denote the degree of vertex $x$ in $G^R$ and $G^B$, respectively.
Thus $d^R_G(x)=|N^R_G(x)|$, $d^B_G(x)=|N^B_G(x)|$. Then $d^R_G(x)+d^B_G(x)=d(x)$. Let $G[S]$ be the subgraph of $G$ induced by $S\subseteq V(G)$.
Note that a subgraph of $G$ admits a red-blue edge coloring preserved from that of $G$.
\medskip

\section{Star-critical Ramsey numbers}
\indent

Slightly abusing notation, for a graph $G$ and a vertex subset $S$ of $G$,
we use $G\setminus S$ for the subgraph induced by $V(G)\setminus S$.

\begin{lemma}\label{l2}
Let $F$ be any bipartite graph on vertex sets $V_1$ and $V_2$, where $s=|V_1|$ and $t=|V_2|$ with $s\le t$,
and $H$ a connected graph of order $n$ that contains a suspended path of length $l$.
If $l\geq (s+1)(s+t-3)+5$, then
$$R_{\mathbb S}(F,H)\geq l-(s+1)(s+t-3)-5.$$
\end{lemma}
{\bf Proof.}~~By Lemma \ref{l1}, we have $r=R(F,H)=n+s-1$.
Let $G$ be graph $K_r\setminus S_N$
with $N=l-(s+1)(s+t-3)-5$.
We shall prove that there is either a red $F$ or a blue $H$ in any red-blue edge coloring of $G$.

Let $v$ be the center vertex of the star $S_N$ which is deleted from $G$ and let $H_1$ be a graph from $H$ with the suspended path shortening by 1.
Since $R(F,H_1)=n+s-2$,
we are done unless there is a blue $H_1$ in $G\setminus \{v\}$.
Let $Y$ be the vertex set $V(G)\setminus V(H_1)$.
Note that $v\in Y$,
and there is a suspended path with $l$ vertices in blue $H_1$, say
$$X'=\{x_1,x_2,\ldots,x_l\}$$
in order. Write $X''=\{x_1,x_2,\ldots,x_{l-1}\}\subseteq X'$.
For any $y\in Y\setminus \{v\}$, consider the $l-1$ edges between $y$ and $X''$.
We assume that no two consecutive edges $yx_i$ and $yx_{i+1}$ are both blue,
since otherwise we have a blue $H$.
Furthermore, suppose that $s+t-1$ edges are blue,
say $yx_{i_1},yx_{i_2},\ldots,yx_{i_{s+t-1}}$ are blue.
Consider any edge $x_{i_{j}+1}x_{i_k+1}$, if it is blue, then there is a blue $H$
with new suspended path
\[
x_1\ldots x_{i_j}yx_{i_k}x_{i_{k}-1}\ldots x_{i_{j}+1}x_{i_k+1}\ldots x_l.
\]
If all edges $x_{i_{j}+1}x_{i_k+1}$ are red, then $x_{i_1+1},x_{i_2+1},\ldots,x_{i_{s+t-1}+1}$
and $y$ will form a red $K_{s+t}$ hence a red $F$.
Consequently, we may assume that any $y\in Y\setminus \{v\}$ is connected with $X''$ in at most $s+t-2$ blue edges.
By the similar argument, vertex $v$ is connected with $X''$ in at most $s+t-1$ blue edges.
As vertex $v$ is adjacent to at least
$(s+1)(s+t-3)+4$
vertices in $X''$,
at least
\[
(s+1)(s+t-3)+4-(s+t-1)-(s-1)(s+t-2)=t
\]
vertices in $X''$ are connected with each vertex of $Y$ in red completely.
These vertices and $Y$ yield a red $F$, completing the proof. \hfill $\square$

\medskip

\begin{lemma}\label{l3}
For graph $F$ of order $m$ with $\chi(F)\geq 3$
and connected graph $H$ containing a suspended path of length $l$,
let $F_1$ be a graph from $F$ by deleting all the $t$ vertices
in a color class among a proper vertex coloring of $F$.
If $l\geq (m-2)(m-t)+t+2$, then
$$R_{\mathbb S}(F,H)\geq \min\{R_{\mathbb S}(F_1,H),l-(m-2)(m-t)-t-2\}.$$
\end{lemma}
{\bf Proof.}~~Set $s=s(F)$, $k=\chi(F)$ and $n=v(H)$.
By Lemma \ref{l1}, we have $r=R(F,H)=(k-1)(n-1)+s$.
Let $G$ be graph $K_r\setminus S_N$ with
$N=\min\{R_{\mathbb S}(F_1,H),l-(m-2)(m-t)-t-2\}$.
We shall prove that there is either a red $F$ or a blue $H$
in any red-blue edge coloring of $G$.
Let $v$ be the center vertex of the star $S_N$ which is deleted from $K_r$
and $H_1$ a graph from $H$ with the suspended path shortening by 1.
Since $R(F,H_1)\le r-1$,
we are done unless there is a blue $H_1$ in $G\setminus \{v\}$.
Delete $n-1$ vertices of this blue $H_1$,
then there are at least $r(F_1,H)$ vertices left and $N\le R_{\mathbb S}(F_1,H)$,
so we may assume that there is a red $F_1$.
Thus we obtain a blue $H_1$ and a red $F_1$.
Let $X$ and $Y$ be their disjoint vertex sets with $|X|=n-1$ and $|Y|=m-t$, respectively.
There is a suspended path with $l$ vertices in blue $H_1$, say
$$X'=\{x_1,x_2,\ldots,x_l\}$$
in order. Write $X''=\{x_1,x_2,\ldots,x_{l-1}\}\subseteq X'$.
Note that $v\notin X'$.
Now we consider two cases.

{\bf Case 1.} If $v\notin Y$,
by the similar argument as in the proof of Lemma \ref{l2},
for any $y\in Y$, $y$ is connected with $X''$ in at most $m-2$ blue edges.
Then there are at most $(m-2)(m-t)$ blue edges between $Y$ and $X''$.
Hence at least
\[
l-1-(m-2)(m-t)>t
\]
vertices in $X''$ are connected with each vertex of $Y$ in red completely,
which yields a red $F$.

{\bf Case 2.} If $v\in Y$,
similarly, for any $y\in Y\setminus \{v\}$,
$y$ is connected with $X''$ in at most $m-2$ blue edges,
and vertex $v$ is connected with $X''$ in at most $m-1$ blue edges.
Since $v$ is adjacent to at least $(m-2)(m-t)+t+1$ vertices in $X''$,
at least
\[
(m-2)(m-t)+t+1-(m-1)-(m-2)(m-t-1)=t
\]
vertices in $X''$ are connected with each vertex of $Y$ in red completely.
These vertices and $Y$ yield a red $F$, completing the proof. \hfill $\square$
\medskip

{\bf Proof of Theorem \ref{t1}.}~~Set $p=v(H)$, $m=v(F)$, $s=s(F)$, $k=\chi(F)$ and $\tau=\tau(F)$.
Let $V_1,V_2,\ldots,V_{k}$ be the color classes of $F$ under a proper vertex coloring using $k$ colors and assume that $|V_1|\ge |V_2|\ge \ldots \ge |V_{k}|$. Choose a vertex coloring such that $|V_k|$ is as small as possible, and then $|V_k|=s$. Let $|V_1|=t$.
We shall show that for sufficiently large $n$,
\begin{equation}\label{e2}
n-p-(m-2)(m-t)-t-1\le R_{\mathbb S}(F,H_n)\leq n+s-\tau-2.
\end{equation}
Then we will have $|R_{\mathbb S}(F,H_n)-n|\le v^2(F)+v(H)$ since $m\ge t+1$.

The upper bound in (\ref{e2}) follows from Lemma \ref{l4} as
$R_{\mathbb S}(F,H_n)\leq n+s-\delta(H)-\tau-1$.
In order to show the lower bound in (\ref{e2}),
as it is trivial for $k=2$,
we assume $k\ge 3$.
Define $F_i=F\setminus (V_1\cup V_2\cup\ldots\cup V_i)$, for $i=1,2,\ldots,k-2$. Note that there is a suspended path in $H_n$ with length $n-p+1$ and $\chi(F_{k-2})=2$.
Applying Lemma \ref{l3}, we have
\[
R_{\mathbb S}(F,H)\ge \min\{R_{\mathbb S}(F_1,H),n-p-(m-2)(m-t)-t-1\}.
\]
Let $|V_{k-1}|=t_1$.
As $m\ge s+t_1+t$,
using Lemma \ref{l3} repeatedly and by Lemma \ref{l2}, we have
\begin{eqnarray*}
R_{\mathbb S}(F,H)&\ge& \min\{R_{\mathbb S}(F_{k-2},H),n-p-(m-2)(m-t)-t-1\}\\
&\ge& \min\{n-p-(s+1)(s+t_1-3)-4,n-p-(m-2)(m-t)-t-1\}\\
&\ge& n-p-(m-2)(m-t)-t-1,
\end{eqnarray*}
which yields the claimed lower bound.\hfill $\square$
\medskip

Now we consider $R_{\mathbb S}(S_m,S_n)$. And the value of $R(S_m,S_n)$ is due to Chv\'atal, Harary \cite{chvatal1972} and Burr, Roberts \cite{burr1973}.
\begin{lemma}\cite{chvatal1972,burr1973}\label{l5}
Let $m$ and $n$ be positive integers. Then
\[
R(S_m,S_n)= \left\{  \begin{array}{cl}
 m+n-1  & \mbox{if $m$ and $n$ are both even,}\\
 m+n & \mbox{otherwise}.
\end{array}      \right.
\]
\end{lemma}

\medskip
{\bf Proof of Theorem \ref{t2}.}~~We consider the following two cases.

\textbf{Case 1.} If $m$ or $n$ is odd. Let $G$ denote an $(S_m,S_n)$-free coloring of $K_{r-1}$ with $r=R(S_m,S_n)=m+n$.
For any vertex $v\in V(G)$, $d^R_G(v)+d^B_G(v)=m+n-2$. Since $G$ is $(S_m,S_n)$-free, for any vertex $v\in V(G)$,
we obtain $d^R_G(v)=m-1$, $d^B_G(v)=n-1$. An additional edge incident with any vertex $v$ in $G$ will produce a red $S_m$ or a blue $S_n$.

\textbf{Case 2.} If $m$ and $n$ are both even. Set $Z_{m+n-2}=\{0,1,\ldots,m+n-3\}$ and $A=\{\pm1,\pm2,\ldots,\pm\frac{m-2}{2}\}$.
Define a graph $H$ on $Z_{m+n-2}$, in which two vertices $x$ and $y$ are adjacent if and only if $x-y\in A$. $H$ is $(m-2)$-regular and its complement $\overline{H}$ is $(n-1)$-regular.
Let $G$ be an edge coloring of $K_{r-1}$ in red and blue with $F^B=\overline{H}\setminus M_{n/2}$ and $r=R(S_m,S_n)=m+n-1$.
Then in graph $G$, there are $m-2$ vertices with red degree $m-2$, blue degree $n-1$,
and $n$ vertices with red degree $m-1$, blue degree $n-2$, respectively. And $G$ is $(S_m,S_n)$-free.
Denote the $m-2$ vertices by set $A_1$ and the $n$ vertices by set $A_2$.
For $S_1=K_2$, an $(S_m,S_n)$-free coloring of $K_{m+n-1}\setminus S_1$ is the graph $G$ and a vertex $v$ adjacent to every vertex in $A_1$ by red edges,
and adjacent to $n-1$ vertices in $A_2$ by blue edges, which implies $R_{\mathbb S}(S_m,S_n)=0$. \hfill $\square$
\medskip

\section{Matching-critical Ramsey numbers}
\indent

{\bf Proof of Theorem \ref{t3}.}~If $m$ and $n$ are both even, then $r=R(S_m,S_n)=m+n-1$.
By Theorem \ref{t2}, we obtain $R_{\mathbb S}(S_m,S_n)=0$ which implies $R_{\mathbb M}(S_m,S_n)=0$.
If $m$ or $n$ is odd, $r=R(S_m,S_n)=m+n$, and then we consider two cases.

\textbf{Case 1.} If $m+n$ is odd, for graph $G=K_r\setminus M_{(m+n-1)/2}$,
then $G$ must contain a vertex, denoted by $v$, such that $d(v)=m+n-1$.
Then we have $d^R_G(v)\geq m$ or $d^B_G(v)\geq n$, producing a red $S_m$ or a blue $S_n$.
Thus $R_{\mathbb M}(S_m,S_n)=(m+n-1)/2$.

\textbf{Case 2.} If $m+n$ is even, for graph $G=K_r\setminus M_{(m+n-2)/2}$,
then $G$ must contain a vertex, denoted by $v$, such that $d(v)=m+n-1$, producing a red $S_m$ or a blue $S_n$.
So we have $R_{\mathbb M}(S_m,S_n)\geq (m+n-2)/2$.
For the upper bound, set $Z_{m+n}=\{0,1,\ldots,m+n-1\}$,
\[
A_1=\Big\{\pm1,\pm2,\ldots,\pm\frac{m-1}{2}\Big\},\;\; A_2=\Big\{\pm\frac{m+1}{2},\pm\frac{m+3}{2},\ldots,\pm\frac{m+n-2}{2}\Big\}.
\]
For graph $H=K_r\setminus M_{(m+n)/2}$, define an $(S_m,S_n)$-free coloring of $H$ on $Z_{m+n}$,
in which the edge $uv$ is red if $u-v\in A_1$ and the edge $uv$ is blue if $u-v\in A_2$.
Then $H^R$ is $(m-1)$-regular and $H^B$ is $(n-1)$-regular, which yields $H$ is $(S_m,S_n)$-free. \hfill $\square$
\medskip

\begin{lemma}\cite{cockayne1,cockayne2,lorimer}\label{l6}
Let $m$ and $n$ be integers with $n\ge m\ge 1$. Then $R(M_m,M_n)=2n+m-1$.
\end{lemma}

\begin{claim}\label{c4}
$K_5\setminus M_2\rightarrow (M_2,M_2).$
\end{claim}
{\bf Proof of Claim \ref{c4}.}~Let $G=K_5\setminus M_2$ with $V(G)=\{v_1,v_2,v_3,v_4,v_5\}$ and $v_1v_2,v_3v_4\notin E(G).$
Suppose that $v_2v_3$ is red, say. Then $v_1v_4$ and $v_1v_5$ are both blue, otherwise there is a red $M_2$. Similarly, we have $v_2v_4$ and $v_3v_5$ are both red, yielding a red $M_2$. \hfill$\square$
\medskip

{\bf Proof of Theorem \ref{t4}.}~We shall prove
$R_{\mathbb M}(M_m,M_n)=\lfloor(2n+m-1)/2\rfloor$ for $n\geq 2$ and $n\geq m\geq 1$. By Lemma \ref{l6}, $r=R(M_m,M_n)=2n+m-1$ for $n\geq m\geq 1$.
We will prove the claimed equality as follows.

\textbf{Case 1.} $m=1$ and $n\geq 2$. In this case, as $r=R(M_1,M_n)=2n$, the claimed equality can be seen easily.

\textbf{Case 2.} $m=n=2$. In this case, by Claim \ref{c4}, any edge coloring of $K_5\setminus M_2$ in red and blue contains either a red $M_2$ or a blue $M_2$.

\textbf{Claim.} If the claimed equality holds for $(M_m,M_n)$, then it holds for $(M_{m+1},M_{n+1})$.

{\bf Proof of Claim.}~Let graph $G=K_{r'}\setminus M_{\lfloor r'/2\rfloor}$ with $r'=R(M_{m+1},M_{n+1})=2n+m+2.$
Note that if every vertex in $G$ has all red edges or all blue edges,
there would be a red $M_{m+1}$ or a blue $M_{n+1}$ in $G$.
Then we may assume that there is a vertex $v$ in $G$ with both a red edge and a blue edge adjacent to vertices $a$ and $b$, respectively.
Let $H=G[V(G)\setminus\{v,a,b\}]$. Since $K_{r}\setminus  M_{\lfloor r/2\rfloor}\subseteq H$ with $r=R(M_{m},M_{n})=2n+m-1$,
$H$ contains either a red copy of $M_{m}$ or a blue copy of $M_{n}$.
Along with edge $va$ or $vb$, we get either a red $M_{m+1}$ or a blue $M_{n+1}$.

Combining Case 1, Case 2 and Claim,
the claimed equality in $(M_m,M_n)$ can be reduced to that in $(M_1,M_{n-m+1})$ or $(M_2,M_2)$, completing the proof. \hfill$\square$
\medskip

\section{Complete-critical Ramsey numbers}
\indent
We prove Theorem \ref{t5} and Theorem \ref{t6} in this section.
\medskip

{\bf Proof of Theorem \ref{t5}.}~Lemma \ref{l5} tells us the Ramsey number $R(S_m,S_n)$.
Assume $m$ or $n$ is odd. As graph $G=K_{m+n}\setminus  K_{m+n-1}=S_{m+n-1}$,
for the center vertex $v\in V(S_{m+n-1})$, we have $d^R_G(v)+d^B_G(v)=m+n-1$.
Then $d^R_G(v)\geq m$ or $d^B_G(v)\geq n$, producing a red $S_m$ or a blue $S_n$.
If $m$ and $n$ are both even, by Theorem \ref{t2}, we have $R_{\mathbb S}(S_m,S_n)=0$, and then we obtain $R_{\mathbb K}(S_m,S_n)=0$. \hfill$\square$
\medskip

{\bf Proof of Theorem \ref{t6}.}~Lemma \ref{l6} yields $r=R(M_m,M_n)=2n+m-1$ for $n\geq m\geq 1$.
For the upper bound, let graph $G=K_{2n+m-1}\setminus K_{n+1}=K_{n+m-2}+(n+1)K_1$
with $G^R=K_{m-1}+2nK_1$ and $G^B=(m-1)K_1\cup \big(K_{n-1}+(n+1)K_1\big)$.
It is easy to see that $G$ is $(M_m,M_n)$-free.

For the lower bound, we will prove the claimed equality as follows.

\textbf{Case 1.} $m=1$ and $n\geq 2$. In this case, any red-blue coloring of the edges of $K_{2n}\setminus K_{n}$ contains either a red $M_1$ or a blue $M_n$.

\textbf{Case 2.} $m=n=2$. In this case, by Claim \ref{c4}, any red-blue coloring of the edges of $K_5\setminus K_2$ contains either a red $M_2$ or a blue $M_2$.

\textbf{Claim.} If the claimed equality holds for $(M_m,M_n)$, then it holds for $(M_{m+1},M_{n+1})$.

{\bf Proof of Claim.}~Let graph $G=K_{2n+m+2}\setminus K_{n+1}=K_{n+m+1}+(n+1)K_1$.
Denote $K_{n+m+1}$ by $G_1$ and $(n+1)K_1$ by $G_2$. Then there must be a vertex $v$ in $G_1$
with two edges in different colors adjacent to vertices $a$ in $G_1$ and $b$ in $G_2$, respectively.
Let $H=G[V(G)\setminus\{v,a,b\}]$. Since $H=K_{r}\setminus  K_{n}$ with $r=R(M_{m},M_{n})=2n+m-1$,
the graph $H$ contains either a red $M_{m}$ or a blue $M_{n}$. Along with edge $va$ or $vb$, we get either a red $M_{m+1}$ or a blue $M_{n+1}$.

Combining Case 1, Case 2 and Claim,
the claimed equality in $(M_m,M_n)$ can be reduced to that in $(M_1,M_{n-m+1})$ or $(M_2,M_2)$, completing the proof. \hfill$\square$

\end{document}